\documentclass[12pt,fleqn]{article}
\usepackage{mathrsfs}
\usepackage{amsfonts}
\usepackage{amssymb}
\usepackage{latexsym,amsmath}
\usepackage{graphicx}
\date{}
\begin{document}
\title{The spectra and the signless Laplacian spectra of graphs with pockets}
\author{Shu-Yu Cui$^a$, Gui-Xian Tian$^b$\footnote{Corresponding author.  E-mail: gxtian@zjnu.cn or guixiantian@gmail.com (G.-X.
Tian)}\\
{\small{\it $^a$ Xingzhi College, Zhejiang Normal University,
Jinhua, Zhejiang, 321004,
P.R. China}}\\
{\small{\it $^b$College of Mathematics, Physics and Information Engineering,}}\\
{\small{\it Zhejiang Normal University, Jinhua, Zhejiang, 321004,
P.R. China}}}\maketitle

\begin{abstract} Let $G[F,V_k,H_v]$ be the graph with $k$ pockets,
where $F$ is a simple graph of order $n\geq1$,
$V_k=\{v_1,\ldots,v_k\}$ is a subset of the vertex set of $F$ and
$H_v$ is a simple graph of order $m\geq2$, $v$ is a specified vertex
of $H_v$. Also let $G[F,E_k,H_{uv}]$ be the graph with $k$
edge-pockets, where $F$ is a simple graph of order $n\geq2$,
$E_k=\{e_1,\ldots,e_k\}$ is a subset of the edge set of $F$ and
$H_{uv}$ is a simple graph of order $m\geq3$, $uv$ is a specified
edge of $H_{uv}$ such that $H_{uv}-u$ is isomorphic to $H_{uv}-v$.
In this paper, we obtain some results describing the signless
Laplacian spectra of $G[F,V_k,H_v]$ and $G[F,E_k,H_{uv}]$ in terms
of the signless Laplacian spectra of $F,H_v$ and $F,H_{uv}$,
respectively. In addition, we also give some results describing the
adjacency spectrum of $G[F,V_k,H_v]$ in terms of the adjacency
spectra of $F,H_v$. Finally, as an application of these results, we
construct infinitely many pairs of signless Laplacian (resp.
adjacency) cospectral graphs.

\emph{AMS classification:} 05C50 05C12

\emph{Keywords:} Adjacency matrix; Signless Laplacian matrix;
Spectrum; Pockets; Edge-pockets

\end{abstract}

\section*{1. Introduction}

Throughout this paper, we consider only finite simple graphs. Let
$G=(V,E)$ be a graph with vertex set $V=\{v_1,v_2,\ldots,v_n\}$ and
edge set $E=\{e_1,e_2,\ldots,e_n\}$. The \emph{adjacency matrix}
$A(G)$ of $G$ is a square matrix of order $n$, whose entry
$a_{i,j}=1$ if $v_i$ and $v_j$ are adjacent in $G$ and $0$
otherwise. Let $D(G)$ be the degree diagonal matrix of $G$. Then the
\emph{Laplacian matrix} $L(G)$ and \emph{signless Laplacian matrix}
$Q(G)$ are defined as $D(G)-A(G)$ and $D(G)+A(G)$, respectively.

For an $n\times n$ matrix $M$ associated to a graph $G$, the
characteristic polynomial $\det(xI_n-M)$ of $M$ is called the
\emph{$M$-characteristic polynomial} of $G$ and is denoted by
$f_M(x)$. The eigenvalues of $M$ (i.e. the zeros of $\det(xI_n-M))$
and the spectrum of $M$ (which consists of the $n$ eigenvalues) are
also called the \emph{$M$-eigenvalues} of $G$ and the
\emph{$M$-spectrum} of $G$, respectively. In particular, if $M$ is
the adjacency matrix $A(G)$ of $G$, then the $A$-spectrum of $G$ is
denoted by $\sigma(A(G)) = (\lambda_1 (G),\lambda_2 (G), \ldots
,\lambda_n (G)), $ where $\lambda_1 (G)\leq \lambda_2
(G)\leq\ldots\leq \lambda_n (G)$ are the eigenvalues of $A(G)$. If
$M$ is the signless Laplacian matrix $Q(G)$ of $G$, then the
$Q$-spectrum of $G$ is denoted by $ \sigma(Q(G)) = (q_1 (G),q_2 (G),
\ldots ,q_n (G)), $ where $q_1 (G)\leq q_2 (G)\leq\ldots\leq q_n
(G)$ are the eigenvalues of $Q(G)$. Throughout this paper, the
$A$-spectrum, $L$-spectrum and $Q$-spectrum denote the adjacency
spectrum, Laplacian spectrum and signless Laplacian spectrum of $G$,
respectively. For more review about the $A$-spectrum, $L$-spectrum
and $Q$-spectrum of $G$, readers may refer to
\cite{Aouchiche2010,Brankov2006,Cvetkovic2007,CvetkovicI,CvetkovicII,CvetkovicIII,Cvetkovic1980,Dam2003,Grone1990,Merris1994,Tian2009}
and the references therein.

The following two definitions come from \cite{Barik2008} and
\cite{Nath2014}, respectively.\\
\\
\textbf{Definition 1.1}\cite{Barik2008} Let $F, H_v$ be graphs of
orders $n$ and $m$, respectively, where $m\geq 2$, $v$ be a
specified vertex of $H_v$ and $V_k=\{v_1,\ldots,v_k\}$ is a subset
of the vertex set of $F$. Let $G=G[F,V_k,H_v]$ be the graph obtained
by taking one copy of $F$ and $k$ vertex disjoint copies of $H_v$,
and then attaching the $i$th copy of $H_v$ to the vertex $u_i$,
$i=1,\ldots,k$, at the vertex $v$ of $H$ (identify $u_i$ with the
vertex $v$ of the $i$th copy). Then the copies of the graph $H_v$
that are attached to the vertices $u_i$, $i=1,\ldots,k$ are referred
to as \emph{pockets},
and $G$ is described as a \emph{graph with $k$ pockets}.\\
\\
\textbf{Definition 1.2}\cite{Nath2014} Let $F$ and $H_{uv}$ be two
graphs of orders $n$ and $m$, respectively, where $n\geq2, m\geq3$,
$E_k=\{e_1,\ldots,e_k\}$ is a subset of the edge set of $F$ and
$H_{uv}$ has a specified edge $uv$ such that $H_{uv}-u$ is
isomorphic to $H_{uv}-v$. Assume that $\mathscr{E}_k$ denote the
subgraph of $F$ induced by $E_k$. Let $G=G[F,E_k,H_{uv}]$ be the
graph obtained by taking one copy of $F$ and $k$ vertex disjoint
copies of $H_{uv}$, and then pasting the edge $uv$ in the $i$th copy
of $H_{uv}$ with the edge $e_i\in E_k$, where $i=1,\ldots,k$. Then
the copies of the graph $H_{uv}$ that are pasted to the edges $e_i$,
$i=1,\ldots,k$ are called as \emph{edge-pockets}, and $G$ is
described as a \emph{graph with $k$ edge-pockets}.\\

Barik\cite{Barik2008} has described the $L$-spectrum of
$G=G[F,V_k,H_v]$ using the $L$-spectra of $F$ and $H_v$, when the
specified vertex $v$ is of degree $m-1$ in $H_v$. In that case, if a
copy of $H_v$ is attached to every vertex of $F$, each at the vertex
$v$ of $H_v$, that is, if $G$ has $n$ pockets, then the graph
$G=G[F,V_n,H_v]$ is nothing but the corona $F\circ H$, where $H =
H_v-v$. Then the complete $L$-spectrum of $G$ is described using the
$L$-spectra of $F$ and $H$\cite{Barik2007}; and if $H$ is a regular
graph or a complete bipartite graph, then the complete $A$-spectrum
and $Q$-spectrum of $G$ are also described using the respective
$A$-spectra and $Q$-spectra of $F$ and
$H$\cite{Barik2007,Cui2012i,Cui2012ii,McLeman2011}.

Recently, Nath and Paul\cite{Nath2014} has described the
$L$-spectrum of $G=G[F,E_k,H_{uv}]$ using the $L$-spectra of $F$ and
$H_v$, when the specified vertices $u$ and $v$ are of degree $m-1$,
and the subgraph $\mathscr{E}_k$ of $F$ induced by $E_k$ is regular.
Similarly, they also describe the $A$-spectrum, when $H_{uv}-\{u,
v\}$ is regular. In that case, if a copy of $H_{uv}$ is pasted to
every edge of $F$, each at the edge $uv$ of $H_{uv}$, that is, if
$G$ has $n$ edge-pockets, then the graph $G=G[F,E_n,H_{uv}]$ is
nothing but the edge-corona $F\diamond H$, where $H =H_{uv}-\{u,
v\}$. Then the complete $L$-spectrum of $G$ is described using the
$L$-spectra of $F$ and $H$ when $F$ is regular\cite{Hou2010}; and if
$F$ is a regular graph and $H$ is also a regular graph or a complete
bipartite graph, then the complete $A$-spectrum and $Q$-spectrum of
$G$ are described using the respective $A$-spectra and $Q$-spectra
of $F$ and $H$\cite{Cui2012i,Cui2012ii,Hou2010}.

Motivated by these researches, we discuss the $Q$-spectrum of
$G=G[F,V_k,H_v]$ and $G[F,E_k,H_{uv}]$. We also consider the
$A$-spectrum of $G=G[F,V_k,H_v]$ when $H = H_v-v$ is regular. The
rest of this paper is organized as follows. In Section 2, we present
some preliminary results, which will be needed to prove our main
results. In Section 3, we give the $A$-characteristic polynomials
and $Q$-characteristic polynomials of $G=G[F,V_k,H_v]$. Using these
results, we describe, except $n + k$ $A$-eigenvalues, all the other
$A$-eigenvalues of $G[F,V_k,H_v]$ in terms of the $A$-eigenvalues of
$H_v$. We also show that the remaining $n+k$ $A$-eigenvalues of
$G[F,V_k,H_v]$ are independent of the graph $H_v$. For the
$Q$-eigenvalues of $G[F,V_k,H_v]$, we also obtain the similar
results. In Section 4, we give the $Q$-characteristic polynomials of
$G[F,E_k,H_{uv}]$ when $\mathscr{E}_k$ be an $r$-regular subgraph of
$F$ induced by $E_k$ in Definition 1.2. Using this result, we
describe, except $n + k$ $Q$-eigenvalues, all the other
$Q$-eigenvalues of $G[F,E_k,H_{uv}]$ in terms of the $Q$-eigenvalues
of $H_{uv}$. We also show that the remaining $n+k$ $Q$-eigenvalues
of $G[F,E_k,H_{uv}]$ are independent of the graph $H_{uv}$. In
addition, we give a complete description of the $Q$-spectrum of
$G[F,E_k,H_{uv}]$ in some particular cases. At the same time, as an
application of these results, we also consider to construct
infinitely many pairs of $A$-cospectral and $Q$-cospectral graphs,
respectively.

\section*{2. Preliminaries}

In this section, we present some preliminary results which will be
needed to prove our main results. In \cite{Cui2012ii}, Cui and Tian
introduced a new invariant, the \emph{$M$-coronal} $\Gamma_M(x)$ of
a matrix M of order $n$ (also see \cite{McLeman2011}). It is defined
to be the sum of the entries of the matrix $(x I_n-M)^{-1}$, that
is,
\[
\Gamma _M (x ) = \textbf{1}_n^T (x I_n - M)^{ - 1} \textbf{1}_n,
\]
where $\textbf{1}_n$ denotes the column vector of size $n$ with all
the entries equal one and $I_n$ is the identity matrix of order $n$.
It is proved\cite{Cui2012ii} that if $M$ is a matrix of order $n$
with each row sum equal to a constant t, then
$\Gamma_M(x)=\frac{n}{x-t}$.

Let $A=(a_{ij})$ and $B=(b_{ij})$ be $m\times n$ and $p\times q$
matrices, respectively. Then the \emph{Kronecker product} of $A$ and
$B$ is defined to the $mp\times nq$ partition matrix $(a_{ij}B)$ and
is denoted by $A\otimes B$. For some properties of the Kronecker
product of matrices, see \cite{Johnson1991}.

Let $G_1$ and $G_2$ be two graphs with disjoint vertex sets $V(G_1)$
,$V(G_2)$ and edge sets $E(G_1)$, $E(G_2)$, respectively. The
\emph{join} $G_1\vee G_2$\cite{Harary1969} of $G_1$ and $G_2$ is the
graph union $G_1\cup G_2$ together with all the edges joining
$V(G_1)$ and
$V(G_2)$.\\
\\
\textbf{Theorem 2.1}\cite{Cvetkovic1980} {\it Let $G_i$ be an $r_i$
regular graph with $n_i$ vertices, where $i=1,2$. Then
\[
f_{A(G_1  \vee G_2 )} (x) = \frac{{f_{A(G_1 )} (x)f_{A(G_2 )}
(x)}}{{(x - r_1 )(x - r_2 )}}\left((x - r_1 )(x - r_2 ) - n_1 n_2
\right).
\]
}\\
\\
\textbf{Theorem 2.2}\cite{Cui2014} {\it Let $G_i$ be an $r_i$
regular graph with $n_i$ vertices, where $i=1,2$. Then
\[
f_{Q(G_1  \vee G_2 )}(x)  = \left( {1 - \frac{{n_1 n_2 }}{{(x - n_1
- 2r_2 )(x  - n_2  - 2r_1 )}}} \right) f_{Q(G_1) } (x - n_2
)f_{Q(G_2) } (x - n_1 ).
\] }\\

Throughout this paper, assume that $F$ is a simple graph of order
$n$, $H_v$ and $H_{uv}$ two simple graphs of order $m$, unless
mentioned otherwise. We also assume the specified vertex $v$ in
$H_v$ is of degree $m-1$  and the specified vertices $u$ and $v$ in
$H_{uv}$ are all of degree $m-1$. Let $H_1 = H_{v}-v$. $H_2 =
H_{uv}-\{u, v\}$. Then $H_v=\{v\}\vee H_1$ and $H_{uv}=(\{u,
v\},\{uv\})\vee H_2$.

If $H_1$ is $r_1$-regular, then by Theorems 2.1 and 2.2, we arrive
at
\begin{equation}\label{1}
\sigma (A(H_v )) = (\sigma (A(H_1 )) - \{ r_1\} ) \cup \{ \alpha
,\beta \},
\end{equation}
where $\alpha ,\beta$ are roots of the equation $x^2-r_1x-m+1=0$,
and
\begin{equation}\label{2}
\sigma (Q(H_v )) = \{ q_j (H_1 ) + 1|1 \le j \le m - 2\}  \cup \{
\gamma,\delta\},
\end{equation}
where $\gamma ,\delta$ are roots of the equation
$(x-2r_1-1)(x-m+1)-m+1=0$.

Similarly, if $H_2$ is $r_2$-regular, then by Theorem 2.2, we obtain
\begin{equation}\label{3}
\sigma (Q(H_{uv} )) = \{ q_j (H_2) + 2|1 \le j \le m - 3\}  \cup
\{m-2, \zeta,\eta\},
\end{equation}
where $\zeta,\eta$ are roots of the equation
$(x-2r_2-2)(x-m)-2(m-2)=0$.

\section*{3. The $A$-spectrum and $Q$-spectrum of $G[F,V_k,H_v]$}

\subsection*{3.1. The $A$-spectrum of $G[F,V_k,H_v]$}

\textbf{Proposition 3.1.} {\it Let $G=G[F,V_k,H_v]$ and $|V_k|=k$.
Then the $A$-characteristic polynomial of $G$ is
\begin{equation}\label{4}
f_{A(G)} (x) = (f_{A(H_1 )} (x))^k \det (xI_n  - M),
\end{equation}
where $M=A(F) + \Gamma _{A(H_1 )} (x)\left( {\begin{array}{*{20}c}
   {I_k } & {0^T }  \\
   0 & 0  \\
\end{array}} \right).
$}\\
\\
\textbf{Proof.} With suitable labeling of the vertices of $G$, we
can write the adjacency matrix of $G$ to
\[
A(G) = \left( {\begin{array}{*{20}c}
   {A(F)} & {\left( {\frac{{I_k  \otimes \textbf{1}_{m - 1}^T }}{0}} \right)}  \\
   {(I_k  \otimes \textbf{1}_{m - 1} |0)} & {I_k  \otimes A(H_1 )}  \\
\end{array}} \right).
\]
Then the $A$-characteristic polynomial of $G$ can be calculated as
follows:
\[
\begin{array}{l}
 f_{A(G)} (x) = \det \left( {\begin{array}{*{20}c}
   {xI_n  - A(F)} & {\left( {\frac{{ - I_k  \otimes \textbf{1}_{m - 1}^T }}{0}} \right)}  \\
   {( - I_k  \otimes \textbf{1}_{m - 1} |0)} & {I_k  \otimes (xI_{m - 1}  - A(H_1 ))}  \\
\end{array}} \right) \\
  \;\;\;\;\;\;\;\;\;\;\;\;\;= \det (xI_{m - 1}  - A(H_1 ))^k \det (S_1 ) \\
 \;\;\;\;\;\;\;\;\;\; \;\;\;= (f_{A(H_1 )} (x))^k \det (S_1 ), \\
 \end{array}
\]
where
\[
\begin{array}{l}
 S_1  = xI_n  - A(F) - \left( {\frac{{I_k  \otimes \textbf{1}_{m - 1}^T }}{0}} \right)\cdot
 \left( {I_k  \otimes (xI_{m - 1}  - A(H_1 ))} \right)^{ - 1} \cdot(I_k  \otimes \textbf{1}_{m - 1} |0) \\
 \;\;\;\; = xI_n  - A(F) - \Gamma _{A(H_1 )} (x)\left( {\begin{array}{*{20}c}
   {I_k } & {0^T }  \\
   0 & 0  \\
\end{array}} \right) \\
 \end{array}
\]
is the Schur complement with respect to $ {I_k  \otimes (xI_{m - 1}
- A(H_1 ))}$. This implies the required result. $\Box$

Let $H_1$ be an $r_1$-regular graph and $|V_k|=k$. Then, except $n +
k$ $A$-eigenvalues, we describe all the other $A$-eigenvalues of
$G[F,V_k,H_v]$ in terms of the $A$-eigenvalues of $H_v$. We also
show that the remaining $n+k$ $A$-eigenvalues of $G[F,V_k,H_v]$ are
independent of the graph $H_v$.\\
\\
\textbf{Theorem 3.2.} {\it Let $H_1$ be an $r_1$-regular graph,
where $r_1\geq1$. Also let
$\lambda\in\sigma(A(H_v))\setminus\{\alpha,\beta\}$, where
$\alpha,\beta$ is described as (\ref{1}), and $G =G[F,V_k,H_v]$. If
$|V_k|=k$, then $\lambda\in\sigma(A(G))$ with multiplicity $k$.
Moreover, the remaining $n + k$ $A$-eigenvalues of
$G$ are independent of $H_v$.}\\
\\
\textbf{Proof.} Since $H_1$ is an $r_1$-regular graph. Thus $\Gamma
_{A(H_1 )}(x)=\frac{m-1}{x-r_1}$ and
\[
f_{A(H_1 )} (x) = (x - r_1)\prod\limits_{j = 1}^{m - 2} {(x -
\lambda _j (H_1 ))}.
\]
Now, from Proposition 3.1, one gets
\begin{equation}\label{5}
f_{A(G)} (x) = (x - r_1)^k \prod\limits_{j = 1}^{m - 2} {(x -
\lambda _j (H_1 ))^k}\det (xI_n  - M),
\end{equation}
where $M=A(F) + \frac{m-1}{x-r_1}\left( {\begin{array}{*{20}c}
   {I_k } & {0^T }  \\
   0 & 0  \\
\end{array}} \right).
$ Notice that $M$ depends on the regularity of $H_1$ only and not on
the structure of $H_1$. From (\ref{1}) and (\ref{5}), we obtain the
required result. $\Box$

As an application of the above results, we may construct many
pairs of $A$-cospectral graphs.\\
\\
\textbf{Corollary 3.3.} {\it Let $H_u$ and $H_v$ be two disjoint
graphs of order $m$ such that $H_u-u$ and $H_v-v$ are $r_1$-regular,
where $m\geq2$ and $r_1\geq1$. If $H_u$ and $H_v$ are
$A$-cospectral, then $G[F,V_k,H_u]$ and $G[F,V_k,H_v]$ are $A$-cospectral.}\\

Let $H_v^*$ be the $H_v$ graph such that $H_1=G_p^{r_1}$, where
$G_p^{r_1}=C_p\Box K_{r_1-1}$ is the Cartesian product of the cycle
$C_p$ and the complete graph $K_{r_1-1}$. Notice that $G_p^{r_1}$ is
$r_1$-regular. Theorem 3.2 implies the following result.\\
\\
\textbf{Corollary 3.4.} {\it Let $H_1$ be an $r_1$-regular graph and
$m-1=p(r_1-1)$ for some integer $p$, where $p\geq3$ and $r_1\geq2$.
Let $G=G[F,V_k,H_v]$ and $G^*=G[F,V_k,H_v^*]$. If $|V_k|=k$, then
$\sigma(A(G))$ consists of the eigenvalues

(a) $\lambda$ with multiplicity $k$, for each
$\lambda\in\sigma(A(H_v))\setminus\{\alpha,\beta\}$, where
$\alpha,\beta$ is described as (\ref{1});

(b)
\[\small
\theta  \in \sigma (A(G^* ))\backslash \{ \underbrace {\lambda _1
(G_p^{r_1} ), \ldots ,\lambda _1 (G_p^{r_1} )}_k,\underbrace
{\lambda _2 (G_p^{r_1} ), \ldots ,\lambda _2 (G_p^{r_1})}_k,\ldots
,\underbrace {\lambda _{m - 2} (G_p^{r_1}), \ldots ,\lambda _{m - 2}
(G_p^{r_1})}_k.
\]}\\
\\
\textbf{Proof.} With suitable labeling of the vertices of $G^*$, we
can write the adjacency matrix of $G^*$ to
\[
A(G^*) = \left( {\begin{array}{*{20}c}
   {A(F)} & {\left( {\frac{{I_k  \otimes \textbf{1}_{m - 1}^T }}{0}} \right)}  \\
   {(I_k  \otimes \textbf{1}_{m - 1} |0)} & {I_k  \otimes A(G_p^r)}  \\
\end{array}} \right).
\]
Thus, from the proofs of Proposition 3.1 and Theorem 3.2, one
obtains
\begin{equation}\label{6}
f_{A(G^*)} (x) = (x - r)^k \prod\limits_{j = 1}^{m - 2} {(x -
\lambda _j (G_p^{r_1}))^k}\det (xI_n  - M),
\end{equation}
where $M=A(F) + \frac{m-1}{x-r_1}\left( {\begin{array}{*{20}c}
   {I_k } & {0^T }  \\
   0 & 0  \\
\end{array}} \right)$. It follows from (\ref{5}) and (\ref{6}) that
\begin{equation}\label{7}
f_{A(G)} (x) = \prod\limits_{j = 1}^{m - 2} {(x - \lambda _j (H_1
))^k } \frac{{f_{A(G^* )} (x)}}{{\prod\limits_{j = 1}^{m - 2} {(x -
\lambda _j (G_p^{r_1}))^k } }}.
\end{equation}
From (\ref{1}) and (\ref{7}), we get $\lambda \in \sigma (A(G))$
with multiplicity $k$, for each
$\lambda\in\sigma(A(H_v))\setminus\{\alpha,\beta\}$. Next we only
need to prove ${\lambda _j (G_p^{r_1})} \in \sigma (A(G^* ))$ with
multiplicity $k$, for each $j=1,\ldots,m-2$.

Let the eigenvalues $\lambda_1(C_p),\ldots,\lambda_p(C_p)=2$ of
$A(C_p)$ are afforded by the eigenvectors
$X_1,\ldots,X_p=\textbf{1}_p$, respectively. Also let the
eigenvalues $\lambda_1(K_{r_1-1}),\ldots$,
$\lambda_{r_1-1}(K_{r_1-1})=r_1-2$ of $A(K_{r_1-1})$ are afforded by
the eigenvectors $Y_1,\ldots,Y_{r_1-1}=\textbf{1}_{r_1-1}$,
respectively. Since $A(G_p^{r_1})=A(C_p)\otimes I_{r_1-1}+I_p\otimes
A(K_{r_1-1})$. Then $X_s\otimes Y_t$ is an eigenvector of
$A(G_p^{r_1} )$ corresponding to the eigenvalue
$\lambda_s(C_p)+\lambda_t(K_{r_1-1})$, where $s=1,\ldots,p$ and
$t=1,\ldots,r_1-1$. Now, for fixed $s$ and $t$, it is easy to verify
that
\[
\left( {\begin{array}{*{20}c}
   0  \\
   {e_1  \otimes X_s  \otimes Y_t }  \\
\end{array}} \right),\left( {\begin{array}{*{20}c}
   0  \\
   {e_2  \otimes X_s  \otimes Y_t }  \\
\end{array}} \right), \ldots ,\left( {\begin{array}{*{20}c}
   0  \\
   {e_k  \otimes X_s  \otimes Y_t }  \\
\end{array}} \right)
\]
are $k$ linearly independent eigenvectors of $A(G^*)$ corresponding
to the eigenvalue $\lambda_s(C_p)+\lambda_t(K_{r_1-1})$, where $ X_s
\otimes Y_t  \ne X_p  \otimes Y_{r_1 - 1} $ and $e_i$ denotes the
column vector of size $k$ with  the $i$-th entry equals one and 0
otherwise. The proof is completed. $\Box$

\subsection*{3.2. The $Q$-spectrum of $G[F,V_k,H_v]$}

\textbf{Proposition 3.5.} {\it Let $G=G[F,V_k,H_v]$ and $|V_k|=k$.
Then the $Q$-characteristic polynomial of $G$ is
\begin{equation}\label{8}
f_{Q(G)} (x) = (f_{Q(H_1 )} (x-1))^k \det (xI_n  - M),
\end{equation}
where $M=Q(F) + (m-1+\Gamma _{Q(H_1 )} (x-1))\left(
{\begin{array}{*{20}c}
   {I_k } & {0^T }  \\
   0 & 0  \\
\end{array}} \right).
$}\\
\\
\textbf{Proof.} With suitable labeling of the vertices of $G$, we
can write the signless Laplacian matrix of $G$ to
\[
Q(G) = \left( {\begin{array}{*{20}c}
   {Q(F) + (m - 1)\left( {\begin{array}{*{20}c}
   {I_k } & {0^T }  \\
   0 & 0  \\
\end{array}} \right)} & {\left( {\frac{{I_k  \otimes \textbf{1}_{m - 1}^T }}{0}} \right)}  \\
   {(I_k  \otimes \textbf{1}_{m - 1} |0)} & {I_k  \otimes (Q(H_1 ) + I_{m - 1} )}  \\
\end{array}} \right).
\]
Then the $Q$-characteristic polynomial of $G$ can be calculated as
follows:
\[
\begin{array}{l}
 f_{Q(G)} (x) = \det \left( {\begin{array}{*{20}c}
   {xI_n  - Q(F) - (m - 1)\left( {\begin{array}{*{20}c}
   {I_k } & {0^T }  \\
   0 & 0  \\
\end{array}} \right)} & {\left( {\frac{{ - I_k  \otimes \textbf{1}_{m - 1}^T }}{0}} \right)}  \\
   {( - I_k  \otimes\textbf{ 1}_{m - 1} |0)} & {I_k  \otimes ((x - 1)I_{m - 1}  - Q(H_1 ))}  \\
\end{array}} \right) \\
 \;\;\;\;\;\;\;\; \;\;\;\;\;\;= \det ((x - 1)I_{m - 1}  - Q(H_1 ))^k \det (S_1 ) \\
  \;\;\;\;\;\;\;\;\;\;\;\;\;\;= (f_{Q(H_1 )} (x - 1))^k \det (S_1 ) \\
 \end{array}
\]

where
\[\small
\begin{array}{l}
 S_1  = xI_n  - Q(F) - (m - 1)\left( {\begin{array}{*{20}c}
   {I_k } & {0^T }  \\
   0 & 0  \\
\end{array}} \right)\\
\;\;\;\;\;\;\;\;- \left( {\frac{{I_k  \otimes \textbf{1}_{m - 1}^T
}}{0}}
\right)(I_k  \otimes ((x - 1)I_{m - 1}  - Q(H_1 )))^{ - 1} (I_k  \otimes\textbf{ 1}_{m - 1} |0) \\
 \;\;\;\;\; = xI_n  - Q(F) - (m - 1)\left( {\begin{array}{*{20}c}
   {I_k } & {0^T }  \\
   0 & 0  \\
\end{array}} \right) - \Gamma _{Q(H_1 )} (x - 1)\left( {\begin{array}{*{20}c}
   {I_k } & {0^T }  \\
   0 & 0  \\
\end{array}} \right) \\
 \end{array}
\]
is the Schur complement with respect to $ {I_k  \otimes ((x - 1)I_{m
- 1}  - Q(H_1 ))}$. Hence, the result follows. $\Box$

Let $H_1$ be a $r_1$-regular graph. If $|V_k|=k$, then except $n +
k$ $Q$-eigenvalues, we describe all the other $Q$-eigenvalues of
$G[F,V_k,H_v]$ using the $Q$-eigenvalues of $H_v$. We also show that
the remaining $n+k$ $Q$-eigenvalues of $G[F,V_k,H_v]$ are
independent of the graph $H_v$.\\
\\
\textbf{Theorem 3.6.} {\it Let $H_1$ be an $r_1$-regular graph,
where $r_1\geq1$. Also let
$q\in\sigma(Q(H_v))\setminus\{\gamma,\delta\}$, where
$\gamma,\delta$ is described as (\ref{2}), and $G =G[F,V_k,H_v]$. If
$|V_k|=k$, then $q\in\sigma(Q(G))$ with multiplicity $k$. Moreover,
the remaining $n + k$ $Q$-eigenvalues of
$G$ are independent of $H_v$.}\\
\\
\textbf{Proof.} If $H_1$ is an $r_1$-regular graph, then $\Gamma
_{Q(H_1 )}(x-1)=\frac{m-1}{x-1-2r_1}$ and
\[
f_{Q(H_1 )} (x-1) = (x -1-2 r_1)\prod\limits_{j = 1}^{m - 2} {(x -1-
q _j (H_1 ))}.
\]
Thus, from Proposition 3.5, one gets
\begin{equation}\label{9}
f_{Q(G)} (x) = (x -1-2 r_1)^k\prod\limits_{j = 1}^{m - 2} {(x -1- q
_j (H_1 ))^k}\det (xI_n  - M),
\end{equation}
where $M=Q(F) + (m-1+\frac{m-1}{x-1-2r_1})\left(
{\begin{array}{*{20}c}
   {I_k } & {0^T }  \\
   0 & 0  \\
\end{array}} \right)
$. Notice that $M$ depends on the regularity of $H_1$ only and not
on the structure of $H_1$. From (\ref{2}) and (\ref{9}), we obtain
the required result. $\Box$

As an application of the above results, we may construct many pairs
of $Q$-cospectral graphs. The following Corollary 3.7 and above
Corollary 3.3 show that there exit two graphs $G[F,V_k,H_u]$ and
$G[F,V_k,H_v]$ such that they are not only $A$-cospectral, but also
are $Q$-cospectral.\\
\\
\textbf{Corollary 3.7.} {\it Let $H_u$ and $H_v$ be two disjoint
graphs of order $m$ such that $H_u-u$ and $H_v-v$ are $r_1$-regular,
where $m\geq2$ and $r_1\geq1$. If $H_u$ and $H_v$ are
$Q$-cospectral, then $G[F,V_k,H_u]$ and $G[F,V_k,H_v]$ are $Q$-cospectral.}\\
\\
\textbf{Corollary 3.8.} {\it Let $H_1$ be an $r_1$-regular graph and
$m-1=p(r_1-1)$ for some integer $p$, where $p\geq3$ and $r_1\geq2$.
Let $G=G[F,V_k,H_v]$ and $G^*=G[F,V_k,H_v^*]$. If $|V_k|=k$, then
$\sigma(Q(G))$ consists of the eigenvalues

(a) $q$ with multiplicity $k$, for each
$q\in\sigma(Q(H_v))\setminus\{\gamma,\delta\}$, where
$\gamma,\delta$ is described as (\ref{2});

(b)

$\theta  \in \sigma (Q(G^* ))\backslash \{ \underbrace {q _1
(G_p^{r_1} )+1, \ldots ,q _1 (G_p^{r_1})+1}_k,\underbrace {q _2
(G_p^{r_1} )+1, \ldots ,q _2 (G_p^{r_1} )+1}_k, \ldots,\\
\indent\indent\indent\underbrace {q _{m - 2} (G_p^{r_1} )+1, \ldots
,q _{m - 2} (G_p^{r_1} )+1}_k\}$.
}\\
\\
\textbf{Proof.} With suitable labeling of the vertices of $G^*$, we
can write the signless Laplacian matrix of $G^*$ to
\[
Q(G^* ) = \left( {\begin{array}{*{20}c}
   {Q(F) + (m - 1)\left( {\begin{array}{*{20}c}
   {I_k } & {0^T }  \\
   0 & 0  \\
\end{array}} \right)} & {\left( {\frac{{I_k  \otimes \textbf{1}_{m - 1}^T }}{0}} \right)}  \\
   {(I_k  \otimes \textbf{1}_{m - 1} |0)} & {I_k  \otimes (Q(G_p^{r_1} ) + I_{m - 1} )}  \\
\end{array}} \right)
\]
Thus, from the proofs of Proposition 3.5 and Theorem 3.6, we obtain
\begin{equation}\label{10}
f_{Q(G^*)} (x) = (x -1-2 {r_1})^k\prod\limits_{j = 1}^{m - 2} {(x
-1- q _j (G_p^{r_1} ))^k}\det (xI_n  - M),
\end{equation}
where $M=Q(F) + (m-1+\frac{m-1}{x-1-2{r_1}})\left(
{\begin{array}{*{20}c}
   {I_k } & {0^T }  \\
   0 & 0  \\
\end{array}} \right)
$. It follows from (\ref{9}) and (\ref{10}) that
\begin{equation}\label{11}
f_{Q(G)} (x) = \prod\limits_{j = 1}^{m - 2} {(x - 1 - q_j (H_1 ))^k
} \frac{{f_{Q(G^* )} (x)}}{{\prod\limits_{j = 1}^{m - 2} {(x - 1 -
q_j (G_p^r ))^k } }}
\end{equation}
From (\ref{2}) and (\ref{11}), we get $q \in \sigma (Q(G))$ with
multiplicity $k$, for each
$q\in\sigma(Q(H_v))\setminus\{\gamma,\delta\}$. Next we only need to
prove ${q _j (G_p^{r_1} )+1} \in \sigma (Q(G^* ))$ with multiplicity
$k$, for each $j=1,\ldots,m-2$.

Let the eigenvalues $q_1(C_p),\ldots,q_p(C_p)=4$ of $Q(C_p)$ are
afforded by the eigenvectors $X_1,\ldots,X_p=\textbf{1}_p$,
respectively, and the eigenvalues $q_1(K_{{r_1}-1}),\ldots$,
$q_{{r_1}-1}(K_{{r_1}-1})=2({r_1}-2)$ of $Q(K_{{r_1}-1})$ are
afforded by the eigenvectors
$Y_1,\ldots,Y_{{r_1}-1}=\textbf{1}_{{r_1}-1}$, respectively. Since
$Q(G_p^{r_1})=Q(C_p)\otimes I_{{r_1}-1}+I_p\otimes Q(K_{{r_1}-1})$.
Then $X_s\otimes Y_t$ is an eigenvector of $Q(G_p^{r_1} )$
corresponding to the eigenvalue $q_s(C_p)+q_t(K_{{r_1}-1})$, where
$s=1,\ldots,p$ and $t=1,\ldots,{r_1}-1$. Now, for fixed $s$ and $t$,
then
\[
\left( {\begin{array}{*{20}c}
   0  \\
   {e_1  \otimes X_s  \otimes Y_t }  \\
\end{array}} \right),\left( {\begin{array}{*{20}c}
   0  \\
   {e_2  \otimes X_s  \otimes Y_t }  \\
\end{array}} \right), \ldots ,\left( {\begin{array}{*{20}c}
   0  \\
   {e_k  \otimes X_s  \otimes Y_t }  \\
\end{array}} \right)
\]
are $k$ linearly independent eigenvectors of $Q(G^*)$ corresponding
to the eigenvalue $q_s(C_p)+q_t(K_{{r_1}-1})+1$, where $ X_s \otimes
Y_t \ne X_p  \otimes Y_{{r_1} - 1} $. The proof is completed. $\Box$

\section*{4. The $Q$-spectrum of $G[F,E_k,H_{uv}]$}

\textbf{Proposition 4.1.} {\it Let $\mathscr{E}_k$ be an $r$-regular
subgraph of $F$ induced by $E_k$ in Definition 1.2. Also let
$G=G[F,E_k,H_{uv}]$ and $|E_k|=k$. Then the $Q$-characteristic
polynomial of $G$ is
\begin{equation}\label{12}
f_{Q(G)} (x) = (f_{Q(H_2 )} (x-2))^k \det (xI_n  - M),
\end{equation}
where $M= Q(F) + r(m - 2)\left( {\begin{array}{*{20}c}
   {I_p } & {0^T }  \\
   0 & 0  \\
\end{array}} \right) + \Gamma _{Q(H_2 )} (x - 2)\left( {\begin{array}{*{20}c}
   {Q(\mathscr{E}_k )} & {0^T }  \\
   0 & 0  \\
\end{array}} \right).
$}\\
\\
\textbf{Proof.} Notice that $\mathscr{E}_k$ has $p=\frac{2k}{r}$
vertices. Let $Q(\mathscr{E}_k)$ be the signless Laplacian matrix of
$\mathscr{E}_k$. With suitable labeling of the vertices of $G$, we
can write the signless Laplacian matrix of $G$ to
\[
Q(G) = \left( {\begin{array}{*{20}c}
   {Q(F) + r(m - 2)\left( {\begin{array}{*{20}c}
   {I_p } & {0^T }  \\
   0 & 0  \\
\end{array}} \right)} & {\left( {\frac{{R(\mathscr{E}_k)  }}{0}} \right)\otimes \textbf{1}_{m - 2}^T }  \\
   {(R(\mathscr{E}_k)^T |0)\otimes \textbf{1}_{m - 2}} & {I_k  \otimes (Q(H_2 ) + 2I_{m - 2} )}  \\
\end{array}} \right).
\]
Then the $Q$-characteristic polynomial of $G$ can be calculated as
follows:
\[
\begin{array}{l}
 f_{Q(G)} (x) = \det \left( {\begin{array}{*{20}c}
   {xI_n  - Q(F) - r(m - 2)\left( {\begin{array}{*{20}c}
   {I_p } & {0^T }  \\
   0 & 0  \\
\end{array}} \right)} & { - \left( {\frac{{R(\mathscr{E}_k )}}{0}} \right) \otimes\textbf{ 1}_{m - 2}^T }  \\
   { - (R(\mathscr{E}_k )^T |0) \otimes \textbf{1}_{m - 2} } & {I_k  \otimes ((x - 2)I_{m - 2}  - Q(H_2 ))}  \\
\end{array}} \right) \\
  \;\;\;\;\;\;\;\;\;\;\;\;\;\;= \det ((x - 2)I_{m - 2}  - Q(H_2 ))^k \det (S_1 ) \\
 \;\;\;\;\;\;\;\;\;\;\;\; \;\;= (f_{Q(H_2 )} (x - 2))^k \det (S_1 ), \\
 \end{array}
\]
where
\[
\begin{array}{l}
 S_1  = xI_n  - Q(F) - r(m - 2)\left( {\begin{array}{*{20}c}
   {I_p } & {0^T }  \\
   0 & 0  \\
\end{array}} \right) \\
 \;\;\;\;\;\;\;- \left( {\frac{{R(\mathscr{E}_k )}}{0}} \right) \otimes \textbf{1}_{m - 2}^T
 \cdot(I_k  \otimes ((x - 2)I_{m - 2}  - Q(H_2 )))^{ - 1}\cdot (R(\mathscr{E}_k )^T |0) \otimes \textbf{1}_{m - 2}  \\
   \;\;\;\;\;= xI_n  - Q(F) - r(m - 2)\left( {\begin{array}{*{20}c}
   {I_p } & {0^T }  \\
   0 & 0  \\
\end{array}} \right) - \Gamma _{Q(H_2 )} (x - 2)\left( {\begin{array}{*{20}c}
   {Q(\mathscr{E}_k )} & {0^T }  \\
   0 & 0  \\
\end{array}} \right) \\
 \end{array}
\]
is the Schur complement with respect to $ {I_k  \otimes ((x - 2)I_{m
- 2}  - Q(H_2 ))}$. This implies the required result. $\Box$

Let $H_2$ be a $r_2$-regular graph. If $|E_k|=k$, then except $n +
k$ $Q$-eigenvalues, we describe all the other $Q$-eigenvalues of
$G[F,E_k,H_{uv}]$ in term of the $Q$-eigenvalues of $H_{uv}$. We
also show that the remaining $n+k$ $Q$-eigenvalues of
$G[F,E_k,H_{uv}]$
are independent of the graph $H_{uv}$.\\
\\
\textbf{Theorem 4.2.} {\it Let $H_2$ be an $r_2$-regular graph,
where $r_2\geq2$. Also let $q\in\sigma(Q(H_{uv}))\setminus\{m-2,
\zeta,\eta\}$, where $\zeta,\eta$ is described as (\ref{3}), and
$G[F,E_k,H_{uv}]$. Assume that $\mathscr{E}_k$ is an $r$-regular
subgraph of $F$ induced by $E_k$. If $|E_k|=k$, then
$q\in\sigma(Q(G))$ with
multiplicity $k$. Moreover, the remaining $n + k$ $Q$-eigenvalues of $G$ are independent of $H_{uv}$.}\\
\\
\textbf{Proof.} Since $H_2$ is an $r_2$-regular graph. Then $\Gamma
_{Q(H_2)}(x-2)=\frac{m-2}{x-2-2r_2}$ and
\[
f_{Q(H_2 )} (x-2) = (x -2-2 r_2)\prod\limits_{j = 1}^{m - 3} {(x -2-
q _j (H_2 ))}.
\]
Thus by Proposition 4.1,
\begin{equation}\label{13}
f_{Q(G)} (x) = (x -2-2 r_2)^k\prod\limits_{j = 1}^{m - 3} {(x -2- q
_j (H_2 ))^k}\det (xI_n  - M),
\end{equation}
where $M= Q(F) + (m - 2)\left( {\begin{array}{*{20}c}
   {rI_p+\frac{Q(\mathscr{E}_k )}{x -2-2 r_2} } & {0^T }  \\
   0 & 0  \\
\end{array}} \right)
$. Notice that $M$ depends on the regularity of $H_2$ only and not
on the structure of $H_2$. From (\ref{3}) and (\ref{13}), we obtain
the required result. $\Box$

As an application of the above results, we may construct many
pairs of $Q$-cospectral graphs.\\
\\
\textbf{Corollary 4.3.} {\it Let $H_{uv}$ and $H_{xy}$ be two
disjoint graphs of order $m$ such that $H_{uv}-\{u,v\}$ and
$H_{xy}-\{x,y\}$ are $r_2$-regular, where $m\geq3$ and $r_2\geq2$.
Let $\mathscr{E}_k$ be an $r$-regular subgraph of $F$ induced by
$E_k$ in Definition 1.2. If $H_{uv}$ and $H_{xy}$ are
$Q$-cospectral, then $G[F,E_k,H_{uv}]$ and $G[F,E_k,H_{xy}]$ are $Q$-cospectral.}\\

Let $H_{uv}^{**}$ be the $H_{uv}$ graph such that $H_2=G_p^{r_2}$,
where $G_p^{r_2}=C_p\Box K_{r_2-1}$ is the Cartesian product of the
cycle $C_p$ and the complete graph $K_{r_2-1}$. Notice that
$G_p^{r_2}$ is
$r_2$-regular. Theorem 4.2 implies the following result.\\
\\
\textbf{Corollary 4.4.} {\it Let $H_2$ be an $r_2$-regular graph and
$m-2=p(r_2-1)$ for some integer $p$, where $p\geq3$ and $r_2\geq2$.
Let $G=G[F,E_k,H_{uv}]$ and $G^{**}=G[F,E_k,H_{uv}^{**}]$. Assume
that $\mathscr{E}_k$ is an $r$-regular subgraph of $F$ induced by
$E_k$. If $|E_k|=k$, then $\sigma(Q(G))$ consists of the eigenvalues

(a) $q$ with multiplicity $k$, for each
$q\in\sigma(Q(H_{uv}))\setminus\{m-2, \zeta,\eta\}$, where
$\zeta,\eta$ is described as (\ref{3});

(b) $\theta  \in \sigma (Q(G^{**} ))\backslash \{ \underbrace {q _1
(G_p^{r_2} )+2, \ldots ,q _1 (G_p^{r_2} )+2}_k,\underbrace {q _2
(G_p^{r_2} )+2, \ldots ,q _2 (G_p^{r_2} )+2}_k,\\
\indent\indent\ldots,\underbrace {q _{m - 3} (G_p^{r_2} )+2, \ldots
,q _{m - 3} (G_p^{r_2} )+2}_k\}$.
}\\
\\
\textbf{Proof.} With suitable labeling of the vertices of $G^{**}$,
we can write the signless Laplacian matrix of $G^{**}$ to
\[
Q(G^{**}) = \left( {\begin{array}{*{20}c}
   {Q(F) + r(m - 2)\left( {\begin{array}{*{20}c}
   {I_p } & {0^T }  \\
   0 & 0  \\
\end{array}} \right)} & {\left( {\frac{{R(\mathscr{E}_k)  }}{0}} \right)\otimes \textbf{1}_{m - 2}^T }  \\
   {(R(\mathscr{E}_k)^T |0)\otimes \textbf{1}_{m - 2}} & {I_k  \otimes (Q(G_p^{r_2}) + 2I_{m - 2} )}  \\
\end{array}} \right).
\]
Thus, from the proofs of Proposition 4.1 and Theorem 4.2, we obtain
\begin{equation}\label{14}
f_{Q(G^{**})} (x) = (x -2-2 r_2)^k\prod\limits_{j = 1}^{m - 3} {(x
-2- q _j (G_p^{r_2}))^k}\det (xI_n  - M),
\end{equation}
where $M= Q(F) + (m - 2)\left( {\begin{array}{*{20}c}
   {rI_p+\frac{Q(\mathscr{E}_k )}{x -2-2 r_2} } & {0^T }  \\
   0 & 0  \\
\end{array}} \right)
$.  It follows from (\ref{13}) and (\ref{14}) that
\begin{equation}\label{15}
f_{Q(G)} (x) = \prod\limits_{j = 1}^{m - 3} {(x - 2- q_j (H_2 ))^k }
\frac{{f_{Q(G^{**} )} (x)}}{{\prod\limits_{j = 1}^{m - 3} {(x - 2 -
q_j (G_p^{r_2} ))^k } }}
\end{equation}
From (\ref{3}) and (\ref{15}), we get $q \in \sigma (Q(G))$ with
multiplicity $k$, for each $q\in\sigma(Q(H_{uv}))\setminus\{m-2,
\zeta,\eta\}$. Next, using the similar technique to the proof of
Corollary 3.8, we may prove that ${q _j (G_p^{r_2} )+2} \in \sigma
(Q(G^{**} ))$ with multiplicity $k$, for each $j=1,\ldots,m-3$. The
proof is completed. $\Box$

In the rest of this paper, we consider many special cases. Assume
that $\mathscr{E}_k$ is an $r$-regular spanning subgraph of $F$.
From the proof of Proposition 4.1, we easily obtain\\
\\
\textbf{Proposition 4.4.} {\it Let $\mathscr{E}_k$ be an $r$-regular
spanning subgraph of $F$. Also let $G=G[F,E_k,H_{uv}]$ and
$|E_k|=k$. Then the $Q$-characteristic polynomial of $G$ is
\begin{equation}\label{16}
f_{Q(G)} (x) = (f_{Q(H_2 )} (x-2))^k \det (xI_n  - M),
\end{equation}
where $M= r((m - 2) + \Gamma _{Q(H_2 )} (x - 2))I_n  + Q(F) + \Gamma
_{Q(H_2 )} (x - 2)A(\mathscr{E}_k )$.}\\

A \emph{$k$-matching} in $G$ is a disjoint union of $k$-edges. If
$2k$ is the order of $G$, then a $k$-matching of $G$ is called a
\emph{perfect matching} of $G$ (\cite{Harary1969}).\\
\\
\textbf{Theorem 4.5.} {\it Let $F=K_{2k}$, $\mathscr{E}_k$ be a
perfect matching of $F$. Also let $H_2$ be an $r_2$-regular graph,
where $r_2\geq2$ and $G=G[F,E_k,H_{uv}]$. Then the $Q$-spectrum of
$G$ is given by

(i) $q$ with multiplicity $k$, for each
$q\in\sigma(Q(H_{uv}))\setminus\{m-2, \zeta,\eta\}$, where
$\zeta,\eta$ is described as (\ref{3});

(ii) $m+2k-4$ with multiplicity $k$;

(iii) two roots of the equation $(x-m-4k+4)(x-2r_2-2)-2(m-2)=0$,
each with multiplicity 1;

(iv) two roots of the equation $(x-m-2k+4)(x-2r_2-2)-2(m-2)=0$, each
with multiplicity $k-1$.
}\\
\\
\textbf{Proof.} From Proposition 4.4, we obtain
\begin{equation}\label{17}
f_{Q(G)} (x) = (x -2-2 r_2)^k\prod\limits_{j = 1}^{m - 3} {(x -2- q
_j (H_2 ))^k}\det (xI_{2k}  - M),
\end{equation}
where
\[
\begin{array}{l}
 M = ((m - 2) + \Gamma _{Q(H_2 )} (x - 2))I_{2k}  + (2k - 2)I_{2k}  + J_{2k}  + \Gamma _{Q(H_2 )} (x - 2)(I_k  \otimes A(K_2 )) \\
  \;\;\;\;\;= \left( {m - 2 + \frac{{m - 2}}{{x - 2 - 2r_2 }} + 2k - 2} \right)I_{2k}  + \frac{{m - 2}}{{x - 2 - 2r_2 }}(I_k  \otimes A(K_2 )) + J_{2k}. \\
 \end{array}
\]
Take $ M_1  = \left( {m - 2 + \frac{{m - 2}}{{x - 2 - 2r_2 }} + 2k -
2} \right)I_{2k}  + \frac{{m - 2}}{{x - 2 - 2r_2 }}(I_k  \otimes
A(K_2 ))$. By a simple computation, one gets
\begin{equation}\label{18}
\begin{array}{l}
 \det (xI_{2k}  - M) = \det (xI_{2k}  - M_1  - J_{2k} ) \\
  \;\;\;\;= \det (xI_{2k}  - M_1 ) \cdot (1 -\Gamma _{M_1 } (x)) \\
  \;\;\;\;= (x - m - 2k + 4)^k \left( {\frac{{(x - m - 2k + 4)(x - 2 - 2r_2 ) - 2(m - 2)}}{{x - 2 - 2r_2 }}} \right)^{k - 1} \left( {\frac{{(x - m - 4k + 4)(x - 2 - 2r_2 ) - 2(m - 2)}}{{x - 2 - 2r_2 }}} \right). \\
 \end{array}
\end{equation}
From (\ref{17}) and (\ref{18}), we obtain the required result.
$\Box$\\
\\
\textbf{Theorem 4.6.} {\it Let $F=K_{n}$, $\mathscr{E}_k=C_n$. Also
let $H_2$ be an $r_2$-regular graph, where $r_2\geq2$ and
$G=G[F,E_k,H_{uv}]$. Then the $Q$-spectrum of $G$ is given by

(i) $q$ with multiplicity $n$, for each
$q\in\sigma(Q(H_{uv}))\setminus\{m-2, \zeta,\eta\}$, where
$\zeta,\eta$ is described as (\ref{3});

(ii) two roots of the equation $(x-2m-2n+6)(x-2r_2-2)-4(m-2)=0$,
each with multiplicity 1;

(iii) two roots of the equation
$(x-2m-n+6)(x-2r_2-2)-2(m-2)(1+\cos\frac{2\pi l}{n})=0$, for each
$l=1,2,\ldots,n-1$.
}\\
\\
\textbf{Proof.} The proof is similar to that of Theorem 4.5,
omitted. $\Box$


{\small \begin{thebibliography}{99}

\bibitem{Aouchiche2010} M. Aouchiche, P. Hansen, A survey of automated
conjectures in spectral graph theory, Linear Algebra Appl., 432
(2010) 2293-2322.

\bibitem{Barik2007} S. Barik, S. Pati, B. K. Sarma, The spectrum of the corona of
two graphs, SIAM. J. Discrete Math., 24 (2007) 47-56.

\bibitem{Barik2008} S. Barik, On the Laplacian spectra of graphs with pockets, Linear
Multilinear Algebra, 56 (2008) 481-490.

\bibitem{Brankov2006} V. Brankov, P. Hansen, D. Stevanovic, Automated upper bounds on
the largest Laplacian eigenvalue, Linear Algebra Appl., 414 (2006)
407-424.

\bibitem{Cui2012i} S-Y. Cui, G-X. Tian, The signless Laplacian spectrum of
the (edge) corona of two graphs, Utilitas Math., 88 (2012) 287-297.

\bibitem{Cui2012ii} S-Y. Cui, G-X. Tian, The spectrum and the signless
Laplacian spectrum of coronae, Linear Algebra Appl., 437 (2012)
1692-1703.

\bibitem{Cui2014} S-Y. Cui, G-X. Tian, The Q-generating function for graphs with
application, arXiv: 1403.2846.

\bibitem{Cvetkovic2007} D. Cvetkovi\'{c}, P. Rowlinson, S.K. Simi\'{c}, Signless Laplacians
of finite graphs, Linear Algebra Appl., 423 (2007) 155-171.

\bibitem{CvetkovicI} D. Cvetkovi\'{c}, S.K. Simi\'{c}, Towards a spectral theory of graphs
based on the signless Laplacian, I, Publ. Inst. Math.(Beogr.), 85
(99) (2009) 19-33.

\bibitem{CvetkovicII} D. Cvetkovi\'{c}, S.K. Simi\'{c}, Towards a spectral theory of graphs
based on the signless Laplacian, II, Linear Algebra Appl., 432
(2010) 2257-2272.

\bibitem{CvetkovicIII} D. Cvetkovi\'{c}, S.K. Simi\'{c}, Towards a spectral theory of graphs
based on the signless Laplacian, III, Applicable Analysis and
Discrete Math., 4 (2010) 156-166.

\bibitem{Cvetkovic1980} D. Cvetkovi\'{c}, M. Doob, H. Sachs, Spectra of
Graphs: Theory and Application, Academic press, New York, 1980.

\bibitem{Dam2003} E. R. Van Dam, W. H. Haemers, Which graphs are
determined by their spectrum?, Linear Algebra Appl., 373 (2003)
241-272.

\bibitem{Grone1990} R. Grone, R. Merris, V S. Sunder, The Laplacian
Spectral of Graphs, SIAM J. Matrix Anal. Appl., 11 (1990) 218-239.

\bibitem{Harary1969} F. Harary, Graph Theory, Addison-Wesley, Reading, MA, 1969.

\bibitem{Hou2010} Y. Hou, W-C. Shiu, The spectrum of the edge corona of two graphs, Electron. J. Linear Algebra., 20 (2010) 586-594.

\bibitem{Johnson1991} R. A. Horn, C. R. Johnson, Topics in matrix analysis, Cambridge University Press, Cambridge, 1991.

\bibitem{Merris1994} R. Merris, Laplacian matrices of graphs: a survey, Linear Algebra
Appl., 197-198 (1994) 143-176.

\bibitem{McLeman2011} C. McLeman, E. McNicholas, Spectra of coronae, Linear Algebra
Appl., 435 (2011) 998-1007.

\bibitem{Nath2014} M. Nath, S. Paul, On the spectra of graphs with edge-pockets, Linear
Multilinear Algebra, DOI: 10.1080/03081087.2013.877010.

\bibitem{Tian2009} G.-X. Tian, T.-Z. Huang, B. Zhou, A note on sum of powers of
the Laplacian eigenvalues of bipartite graphs, Linear Algebra Appl.,
430 (2009) 2503-2510.

\end{thebibliography}}
\end{document}